\documentclass[11pt]{article}
\usepackage{amssymb}
\usepackage{amsmath}

\newtheorem{theo}{Theorem}
\newtheorem{prop}{Proposition}
\newtheorem{lemm}{Lemma}
\newtheorem{coro}{Corollary}
\newtheorem{rema}{Remark}

\newtheorem{ex}{Example}

\newcommand{\cqfd}
{%
\mbox{}%
\nolinebreak%
\hfill%
\rule{2mm}{2mm}%
\medbreak%
\par%
}
\newfont{\gothic}{eufb10}
\date{\empty}
\begin{document}
\title{A generalization of the Kuga-Satake construction }
\author{Claire Voisin\\ Institut de math\'{e}matiques de Jussieu, CNRS,UMR
7586} \maketitle \setcounter{section}{-1}
\begin{flushright} {\it \`A la m\'emoire d'Armand Borel}
\end{flushright}
\section{Introduction}
The Kuga-Satake construction \cite{kugasatake} associates to a
polarized Hodge structure $H$ of weight $2$ with $h^{2,0}=1$ an
abelian variety $A$ which satisfies the property that $H$ is a
sub-Hodge structure of $Hom\,(H_1(A),H_1(A))$. The construction is
very tricky and intriguing geometrically: one first associates to
the lattice $(H,<,>)$ its Clifford algebra $C(H)$, which is again a
lattice. Then one constructs a complex structure on
$C(H)\otimes\mathbb{R}$, using the rank $1$ subspace $H^{2,0}\subset
H\otimes\mathbb{C}$ defining the Hodge structure on $H$. Thus the
quotient
$$\frac{C(H)\otimes\mathbb{R}}{C(H)}$$ is endowed with the structure of a complex torus, and with some more work, on
can show that it is in fact an abelian variety. This abelian variety
$A$ has by definition $H_1(A,\mathbb{Z})=C(H)$ and the morphism of
weight $2$ Hodge structures
$$H\rightarrow End\,(H_1(A,\mathbb{Z}))$$
is given by Clifford multiplication on the left acting on $C(H)$: $H\rightarrow End\,(C(H))$.

In \cite{deligne}, Deligne proved that a general weight $2$ polarized  Hodge structure coming from
geometry, is not a quotient of
a Hodge structure of the form $K\otimes L$, where $K$ and $L$ are weight $1$
polarized Hodge structures.
His argument is  that the Mumford-Tate group \cite{vg} of a Hodge structure
$K\otimes L$ has a very restricted form, while the Mumford-Tate group
for general weight $2$  Hodge structure coming from
geometry is very large, as it contains the monodromy group, and thus
can be in some cases the whole orthogonal group
(cf \cite{beauville}).

In the papers \cite{paranjape}, \cite{voisinsym}, \cite{vangeemen}, some
constructions are given, which realize geometrically for $K3$-surfaces with large Picard number
 the Kuga-Satake construction, as expected from the Hodge conjecture:
the inclusion
$$ H\subset Hom\,(H_1(A),H_1(A))\cong H^1(A)\otimes H^1(A)$$
of weight $2$ rational Hodge structures can be understood as a degree $4$ Hodge class
in $H^2(S)\otimes H^2(A\times A)$, and thus should correspond to a codimension $2$ cycle
$$Z\subset S\times A\times A$$
with rational coefficients, such that the inclusion above is given
by
$$[Z]_*: H^2(S)\rightarrow H^2(A\times A).$$
 Morrison \cite{morrison} proves that if $S$ is a Kummer surface,
 that is the minimal desingularization of the quotient of an abelian surface by the $-Id$ involution:
$$S=\widetilde{T/\pm1},$$ then its Kuga-Satake variety
$A$ is a sum of copies of $T$.
Paranjape \cite{paranjape} solves the problem for the members of a certain  family of $K3$-surfaces with
Picard number $16$, by proving that such a $K3$ surface $S$ is dominated by
the self-product of curve $C$, such that the Kuga-Satake variety of $S$
is a sum of copies of an abelian subvariety of $J(C)$.

Our goal in this paper is rather to generalize the Kuga-Satake construction to other
weight $2$ Hodge structures. It is based on  the following observation :
There are in fact two steps in the Kuga-Satake construction: The first one is the Clifford
algebra construction,
and it turns out that
if a polarized weight $2$ Hodge structure $H$ has $h^{2,0}=1$, then its Clifford
algebra $C(H)$ also carries a polarized weight $2$ Hodge structure (see  section \ref{subsec2}), which is
compatible with the algebra structure.
 The second step consists in associating to
 the polarized weight $2$ Hodge structure, which is
compatible with the Clifford algebra structure, a weight $1$ Hodge structure.
We generalize here this second step as follows:
\begin{theo}\label{trivialcenter} Let
$H$ be a non trivial polarized integral Hodge structure of weight $2$. Assume there is
a (associative, unitary) ring structure on $H$, satisfying the following conditions:
\begin{enumerate}
\item \label{item1}The product $H\otimes H\rightarrow H$ is a morphism of Hodge structures of bidegree $(-1,-1)$.
\item \label{item2} There is an (involutive) endomorphism   $t:H\rightarrow H$ such that
for any $a\in H$, multiplication (on left and right) by $t(a)$ is  adjoint to multiplication
(on left and right) by $a$ with respect to the intersection form $<,>$ which gives the polarization.
\item \label{item3} The $t$-invariant part of the center of $H$ is $\mathbb{Z}$.
\end{enumerate}
Then there exists a unique weight $1$
 Hodge structure on $H$ which is polarizable, (or equivalently an  abelian variety
 $A$ with $H=H_1(A,\mathbb{Z})$,) such that
multiplication on the left
$$H\rightarrow Hom\,(H,H)$$
is a morphism of Hodge structures, where on the right, we put the weight $2$ Hodge structure
 induced by the weight
$1$ Hodge structure on $H$.
\end{theo}
We explain this construction in section \ref{sec1}. In section
\ref{sec2}, we  show that under the assumptions above,
any simple abelian variety $B$ such that there is
an inclusion  of weight $2$ Hodge structures, which is also the inclusion of a subring
stable under some Rosatti involution,
$$H\subset Hom\,(H_1(B),H_1(B))$$
must be a quotient of the abelian variety $A$ associated to $H$ by Theorem
\ref{trivialcenter}. We also compare this construction
with the Kuga-Satake construction.

In the final section, we turn to the problem of removing the assumption on the center.
Passing to rational coefficients,
we analyse the structure of the center $K$ of an algebra satisfying
properties \ref{item1}, \ref{item2} above. We show that it is a product of
number fields $K_i$ which are either  totally real fields or a quadratic extension
of  totally real fields, and that unless some corresponding factor
$H_i$ of $H$ is a  simple central algebra over $K_i$  the same conclusion as in Theorem
\ref{trivialcenter} holds.
We finally study the  last case. We then show that
we always have existence of an abelian variety $A$ as in Theorem \ref{trivialcenter} but not uniqueness.

\vspace{0,5cm}

{\bf Thanks.} This work was completed at ETH,  which
I wish to thank  for its hospitality. I also would like  to thank the organizers of the Armand Borel
memorial conference in Hangzhou for allowing me to present these results, and especially
Lizhen Ji for his encouragements to write  them up for this  volume in memory of A. Borel.

\section{The case of trivial center \label{sec1}}
 We start with a polarized integral Hodge structure $H$ of weight $2$, with
 $H^{2,0}\not=0$. Thus
 $H$ is a lattice, endowed with a decomposition
 $$H_\mathbb{C}:=H\otimes\mathbb{C}=H^{2,0}\oplus H^{1,1}\oplus H^{0,2},$$
 $$H^{0,2}=\overline{H^{2,0}},\,H^{1,1}=\overline{H^{1,1}},$$
 and an integral symmetric bilinear form
 $<,>$ satisfying the following properties
 (Hodge-Riemann bilinear relations):

 i) The Hodge decomposition is orthogonal with respect to the Hermitian intersection pairing
 $h(\alpha,\beta)=<\alpha,\overline{\beta}>$ on $H\otimes\mathbb{C}$.

 ii) The Hermitian form $h$ is positive definite  on $H^{2,0}$ and $H^{0,2}$, negative definite on
 $H^{1,1}$.

 We assume furthermore that
 $H$ has the structure of an unitary associative ring, and that conditions
 \ref{item1}, \ref{item2} of Theorem \ref{trivialcenter} are satisfied.
For
 $K,\,L$ two vector subspaces of $H_\mathbb{C}$, $KL$ will be the subspace generated by
 products $kl$ of elements of $K$ and $L$.
Condition \ref{item1} means concretely that
$$H^{2,0}H^{2,0}=0,\, H^{2,0}H^{1,1}\subset H^{2,0}, \,H^{1,1}H^{2,0}\subset H
^{2,0},\,H^{1,1}H^{1,1}\subset H^{1,1},$$
$$H^{0,2}H^{0,2}=0,\, H^{1,1}H^{0,2}\subset H^{0,2},\,
H^{0,2}H^{1,1}\subset H^{0,2}.$$
Condition \ref{item2} means that the intersection form
$<,>$ can be written as
$$<a,b>=<at(b),1>=<t(b)a,1>,\,\forall a,\,b\in H.$$

We observe first of all that the map $t$ (which is a
morphism of rings $H\rightarrow H^{op}$) has to be also a morphism of Hodge structures, that
is
$$t(H^{p,q})\subset H^{p,q}.$$
Indeed, this follows immediately from the fact that the product map and the polarization are compatible with the
Hodge decomposition.

 We consider the subspace
 $$W=H^{2,0}H_\mathbb{C}\subset H_\mathbb{C}.$$
 We have the following Proposition.
 \begin{prop} i)  We have $W\cap\overline{W}=0$.

 ii) The sum $W\oplus\overline{W}$ is a two-sided ideal of $H_\mathbb{C}$ which is stable
 under the Hodge decomposition and under $t$.

 \end{prop}
{\bf Proof.} i) Note first that $W\subset F^1H:=H^{2,0}\oplus H^{1,1}$ and similarly
$\overline{W}\subset H^{1,1}\oplus H^{0,2}$. Thus
$$W\cap\overline{W}\subset H^{1,1}.$$

Let $x\in W\cap\overline{W}$. By definition, this $x$ can be written as
$$x=\sum_i\alpha_i\beta_i,$$ with $\alpha_i\in H^{2,0}$ and $\beta_i\in H^{0,2}$, and as
$$x=\sum_j\gamma_j\delta_j,$$ with $\gamma_j\in H^{0,2}$ and $\delta_j\in H^{2,0}$.
Applying $t(\overline{\cdot})$ to the second expression, we get
$$t(\overline{x})=\sum_jt(\overline{\delta_j})t(\overline{\gamma_j}),$$
with $t(\overline{\delta_j})\in H^{0,2},\,t(\overline{\gamma_j})\in H^{2,0}$.
Thus, as $H^{2,0}H^{2,0}=0$, we conclude that
$$t(\overline{x})x=0.$$
It follows that
$$<t(\overline{x})x,1>=0=<x,\overline{x}>,$$
and as $x\in H^{1,1}$, this implies $x=0$ by the condition ii) above satisfied by the polarization.

\vspace{0,5cm}

ii) By definition, $W$ and $\overline{W}$ are stable under right multiplication.
Next, if $a\in H$ and $b=\eta w\in W$, with $\eta\in H^{2,0}$,
then we can write
$a=a_1+a_2$, where
$$a_1\in F^1H_\mathbb{C},\,a_2\in H^{0,2}.$$
We then have
$$ ab=a_1\eta w+ a_2b,$$
and as $F^1H_\mathbb{C} H^{2,0}\subset H^{2,0}$,
we have $a_1\eta\in H^{2,0}$. Thus the first term belongs to $W$ and the second belongs to
$\overline{W}=H^{0,2}H_\mathbb{C}$. This shows that
$W\oplus\overline{W}$ is a two-sided ideal.

Finally, note that, because $H$ is unitary, $W$ contains $H^{2,0}$, and thus can be written as
$$ W=H^{2,0}\oplus H^{2,0}H^{0,2},$$
because the second space is contained in $H^{1,1}$. This shows that $W$ and hence also
$W\oplus\overline{W}$, are stable under Hodge decomposition.
As the space $H^{2,0}$ is stable under $t$,
we have
$$t(H^{2,0}H^{0,2})=H^{0,2}H^{2,0}=\overline{H^{2,0}}\,\overline{H^{0,2}}\subset \overline{W},$$
hence we get $t(W)\subset W\oplus\overline{W}$, and thus
$t(W\oplus\overline{W})\subset W\oplus\overline{W}$.

\cqfd

\begin{coro} We have a decomposition of the $\mathbb{C}$-algebra
$H_\mathbb{C}$ as a direct sum
$$H_\mathbb{C}=(W\oplus\overline{W})\oplus M,$$
where $M$ is defined as the orthogonal
of $W\oplus\overline{W}$ with respect to $<,>$.
\end{coro}
{\bf Proof.} Indeed, we know that $W\oplus \overline{W}$ is stable under the Hodge decomposition and under complex
conjugation. It follows that the intersection form
$<,>$ is non degenerate on $W\oplus\overline{W}$, because the Hermitian form
$h$ is non degenerate on each of its $(p,q)$-piece and the Hodge decomposition is orthogonal for $h$.
Thus we have a orthogonal decomposition
\begin{eqnarray}\label{orthog} H_\mathbb{C}=(W\oplus\overline{W})\oplus M.
\end{eqnarray}
As $W\oplus\overline{W}$ is a two-sided ideal which is stable under
$t$, it follows that the same is true for its orthogonal $M$.
Thus we must have $wm=0$ for $w\in (W\oplus \overline{W}),\,m\in M$, which shows that
(\ref{orthog}) is a decomposition of the algebra as a direct sum.
\cqfd
\begin{coro} \label{presque} If the $t$-invariant part of the center of  $H$  is
equal to $\mathbb{Z}$, then
$M=0$. Thus
\begin{eqnarray}\label{hodgedecomp}H_\mathbb{C}=W\oplus\overline{W}
\end{eqnarray}
which defines a complex structure on $H_\mathbb{R}$.

\end{coro}
{\bf Proof.} Indeed, under this assumption, we also get that the
$t$-invariant part of the center of the $\mathbb{C}$-algebra
$H_\mathbb{C}$ is $\mathbb{C}$. But as the decomposition
(\ref{orthog}) is orthogonal and is a decomposition into $t$-invariant subspaces, the associated idempotents
$e_{W\oplus\overline{W}}$ and $e_M$ are $t$-invariant and of course they are central.
Thus, under the assumption of  Corollary \ref{presque}, we have either
$e_{W\oplus\overline{W}}=0$ or $e_M=0$. As $H^{2,0}\not=0$, we have $W\not=0$ and thus
$e_{W\oplus\overline{W}}\not=0$, hence $e_M=0$ and $M=0$.

\cqfd
Corollary \ref{presque} proves part of the existence statement of Theorem \ref{trivialcenter}. In order to have
the complete proof of existence, we need to show that the Hodge structure of weight $1$ defined
by the decomposition \ref{hodgedecomp} is polarizable (or equivalently, that the
corresponding complex torus
$$A=H_\mathbb{C}/(W\oplus H)$$ is an abelian variety).
 This is done as follows:
Let $a\in H$ be such that $t(a)=-a$.
Consider the skew pairing
$$\omega_a(w,w')=<w,w'a>$$
on $H$ and denote in the same way
its $\mathbb{C}$-linear extension to $H_\mathbb{C}$.
\begin{lemm} The subspace $W\subset H_\mathbb{C}$ is totally isotropic with respect to
$\omega_a$.
\end{lemm}
{\bf Proof.} let $w=\eta m,\,w'=\eta'm'\in W$, with $\eta\in H^{2,0},\,\eta'\in H^{2,0}$.
Then
$$\omega_a(w,w')=<w,w'a>=<t(a)t(w')w,1>=<t(a)t(m')t(\eta')\eta m,1>,$$
and this is $0$ because $t(\eta')\eta=0$.
\cqfd
Let us now show the following
\begin{prop}\label{polar} For an adequate choice of $a$, $\omega_a$ polarizes the weight $1$
Hodge structure given by the decomposition (\ref{hodgedecomp}).
\end{prop}
{\bf Proof.}
In order that  $\omega_a$ defines a polarization, it needs to satisfy the property that
the Hermitian form $h_a$ defined by
$$ h_a(w,w')=i\omega_a(w,\overline{w'})$$
is  positive definite on $W$.
As rational elements $a\in H_\mathbb{Q}$ are dense in $H_\mathbb{R}$, it suffices to show that for some
$a\in H_\mathbb{R}$, this property is satisfied.
Let us take $a$ to be a sum of  terms  of the following form:
$$a^-= -i(\eta\overline{\eta}-\overline{\eta}\eta),
$$
where $\eta\in H^{2,0}$ satisfies $t(\eta)=-\eta$ and
$$a^+=i(\eta\overline{\eta}-\overline{\eta}\eta),
$$
where $\eta\in H^{2,0}$ satisfies $t(\eta)=\eta$. Then
$a$ is real and $t(a)=-a$.

We have to compute the
sign of the Hermitian form $h_a$ on
$W=H^{2,0}\oplus H^{2,0}H^{0,2}$. Note that
for $\mu\in H^{2,0}$ and $\nu\in H^{2,0}H^{0,2}$,
we have
$$h_{a^+}(\mu,\nu)=i<\mu,\overline{\nu}a^+>=
-<\mu,\overline{\nu}(\eta\overline{\eta}-\overline{\eta}\eta)>$$
$$=<\mu,\overline{\nu}\,\overline{\eta}\eta>$$
because $\overline{\nu}\eta=0,\,\forall \eta\in H^{2,0}$. The last term is equal to
$$<\mu t(\eta),\overline{\nu}\,\overline{\eta}>$$
which is also $0$
because $\mu t(\eta)=0,\,\forall \eta\in H^{2,0}$.

 The same computation works with $
h_{a^-}$ and it thus follows that
$H^{2,0}$ and $H^{2,0}H^{0,2}$ are perpendicular with respect
to $h_a$.  Hence it suffices to compute the signs of $h_{a}$ on each term
$H^{2,0}$ and $H^{2,0}H^{0,2}$.

If now $\mu\in H^{2,0}$, we have
$$h_{a^{\pm}}(\mu,\mu)=i<\mu,\overline{\mu}a^
\pm>=
-\pm<\mu,\overline{\mu}(\eta\overline{\eta}-\overline{\eta}\eta) >$$
$$=-\pm<\mu,\overline{\mu}\eta\overline{\eta}>,
$$
because $\overline{\mu}\,\overline{\eta}=0$, for all $\eta\in H^{2,0}$, and this is equal to
\begin{eqnarray}
\label{term1}-\pm<\mu t(\overline{\eta}),\overline{\mu}\eta>.
\end{eqnarray}
In the case of $a^+$ where $t(\eta)=\eta$, and in the case of
$a^-$ where $t(\eta)=-\eta$, we get
$$-<\mu \overline{\eta},\overline{\mu \overline{\eta}}>.$$
But $\mu \overline{\eta}\in H^{1,1}$ and the Hermitian form
$<\alpha,\overline{\beta}>$ is negative on $H^{1,1}$.
Thus (\ref{term1}) is $\geq0$.

Let us now consider the case of
$\nu\in H^{2,0}H^{0,2}$. In this case we have
$$h_{a^{\pm}}(\nu,\nu)=i<\nu,\overline{\nu}a^\pm>=
-\pm<\nu,\overline{\nu}(\eta\overline{\eta}-\overline{\eta}\eta) >$$
$$=\pm<\nu,\overline{\nu}\,\overline{\eta}\eta >,$$
because $\overline{\nu}\eta=0,\,\forall \eta\in H^{2,0}$.
This is also equal to
\begin{eqnarray}
\label{term2}\pm<\nu t(\eta),\overline{\nu}\,\overline{\eta} >.
\end{eqnarray}
In the case of $a^+$ where $t(\eta)=\eta$, and in the case of
$a^-$ where $t(\eta)=-\eta$, we get
$$<\nu \eta,\overline{\nu\eta} >.$$
But $\nu\eta\in H^{2,0}$ and the Hermitian form
$<\alpha,\overline{\beta}>$ is positive on $H^{2,0}$. Thus
(\ref{term2}) is also $\geq0$.

To conclude, it remains to show that for $a$ a generic sum of terms $a^+$,
$a^-$ as above, $h_a$ is non degenerate, or equivalently does not vanish
on $W\setminus \{0\}$.
But the computation above shows that a null-vector $\mu$ of
a generic sum $\sum h_{a^{\pm}}$ in $H^{2,0}$ has to satisfy
$$\mu \overline{\eta}=0,\,\forall \eta\in H^{0,2}$$
and that a
null-vector $\nu$ of a generic sum of
$h_{a^{\pm}}$ in $H^{2,0}H^{0,2}$ has to satisfy
$$\nu \eta=0,\,\forall \eta\in H^{2,0}.$$
In the first case, we get $\mu t(\overline{\mu})=0$, which implies that
$$<\mu t(\overline{\mu}),1>=<\mu,\overline{\mu}>=0.$$
As $\mu\in H^{2,0}$, this implies by the second Hodge-Riemann bilinear relations
that $\mu=0$.

In the second case, as $t(\overline{\nu})\in H^{2,0}H^{0,2}$, we conclude that
$\nu t(\overline{\nu})=0$, which implies that
$$<\nu t(\overline{\nu}),1>=<\nu,\overline{\nu}>=0.$$
As $\nu\in H^{1,1}$, this implies by the second Hodge-Riemann bilinear relations
that $\nu=0$.

\cqfd

In order to conclude the proof of Theorem \ref{trivialcenter}, it suffices to prove the uniqueness
statement. But this is clear, because if
a weight $1$ Hodge structure on $H$ given by a decomposition
\begin{eqnarray}\label{eqn}H_\mathbb{C}=W'\oplus\overline{W'}
\end{eqnarray}
satisfies the property that left multiplication
$$H\rightarrow Hom\,(H,H)$$
is a morphism of weight $2$ Hodge structures, then
we must have
$$H^{2,0}H_\mathbb{C}\subset W'.$$
Indeed, the
$(2,0)$-piece of $Hom\,(H,H)$ for the weight $2$ Hodge structure
induced by the weight $1$ Hodge structure on $H$ given by (\ref{eqn}) is equal to
$Hom\,(\overline{W'},W')$ and thus is contained
in $Hom_\mathbb{C}\,(H_\mathbb{C},W')$.

Thus by definition of
$W$, we must have $W\subset W'$, and then equality for dimension reasons.

\cqfd

\section{General properties, examples \label{sec2}}
\subsection{General properties}
We first start with the proof of the following:
\begin{prop}\label{universal} Let $B$ be a simple abelian variety, and let
$$H\subset End\,(H_1(B,\mathbb{Z}))$$
be a subring which is also a sub-Hodge structure, satisfying the conditions
of Theorem \ref{trivialcenter}.
Then $B$ is  a quotient of the abelian variety $A$
constructed in Theorem \ref{trivialcenter}.
\end{prop}
(Actually, only conditions \ref{item2} and \ref{item3}
are to be verified. \ref{item2} means that $H$ has to be stable under one Rosatti involution.)

\vspace{0,5cm}

{\bf Proof.} We want to show equivalently that there exists a non trivial morphism of Hodge
structure
$$\alpha:H\rightarrow H_1(B,\mathbb{Z})$$
where on the left, $H$ is endowed with the weight $1$ Hodge structure given in
Theorem \ref{trivialcenter} and on the right, the Hodge decomposition
is dual to the Hodge decomposition on $H^1(B)$ (or equivalently is the Hodge decomposition on
$H^{2n-1}(B)$, $n=dim\,B$) :
$$H_1(B,\mathbb{C})=H_{1,0}(B)\oplus H_{0,1}(B),$$
with
$$B=H_1(B,\mathbb{C})/(H_{1,0}(B)\oplus H_1(B,\mathbb{Z})).$$
 Indeed, as $B$ is simple, the induced non trivial
morphism of abelian varieties
$$\alpha:A\rightarrow B$$
 has to be surjective.

Let $\beta\in H_1(B,\mathbb{Z})$ and consider the
map
$$e_\beta: H\rightarrow H_1(B,\mathbb{Z}),$$
$$h\mapsto h(\beta).$$
Certainly this map is non zero for at least one $\beta$.

We claim that this is a morphism of Hodge structures. Indeed, we only have to show that
$e_\beta(W)\subset H_{1,0}(B)$. But, as $H$ is a sub-Hodge structure
of $End\,H_1(B)$,
$$H^{2,0}\subset H^{2,0}(End_\mathbb{C}\,(H_1(B,\mathbb{C}))
\subset Hom_\mathbb{C}\,(H_1(B,\mathbb{C}),H_{1,0}(B)).$$
Thus $W=H^{2,0}H_\mathbb{C}\subset Hom_\mathbb{C}\,(H_1(B,\mathbb{C}),H_{1,0}(B))$, which
proves the claim.

\cqfd
\begin{rema}{\rm  Without the assumption that
$B$ is simple, we still  get the following statement :
if the map $H\otimes H_1(B,\mathbb{Q})\rightarrow H_1(B,\mathbb{Q})$,
$h\otimes \beta\mapsto h(\beta)$, is surjective, then
$B$ is a quotient of a sum of copies of $A$. }

\end{rema}
Indeed, choosing a basis $\beta_i,\,i=1,\ldots,n$ of $H_1(B,\mathbb{Z})$, the
morphism
$$\sum_ie_{\beta_i}:A^n\rightarrow B$$
is then surjective.
\cqfd
Finally, observe that the argument above can be reversed to show the following:
\begin{prop}Let $B$ be an abelian variety, and let
$$H\subset End\,(H_1(B,\mathbb{Z}))$$
be a subring which is also a sub-Hodge structure, satisfying the conditions
of Theorem \ref{trivialcenter}.
Then  the associated abelian variety $A$
is isogenous to an abelian subvariety
of a sum of copies of $B$.

\end{prop}

Indeed, we use the morphisms of Hodge structures $e_\beta$, which give morphism of abelian varieties
$A\rightarrow B$. As $H\subset End\,(H^1(B,\mathbb{Z}))$, for $\beta_i,\,i=1,\ldots,n$, running over a basis
of $H_1(B,\mathbb{Z})$, this gives a morphism of abelian varieties
$$(e_{\beta_i}):A\rightarrow B^n$$
 which has a finite kernel.
\cqfd
\subsection{The Kuga-Satake construction.\label{subsec2}}

\vspace{0,5cm}

We start from a polarized Hodge structure $(H,<,>)$ of weight $2$.
Consider the Clifford algebra $C(H)$ which is the quotient of the
tensor algebra $\bigotimes H$ by the ideal generated by the
relations
$$h\otimes h=-<h,h>1.$$
For sign reasons, it is better to work with the even part $C^+(H)$
generated by products of an even number of elements of $H$, but we
won't do this, as it makes computations more complicated. As a
lattice, $C(H)$ is canonically isomorphic to the  exterior algebra
$\bigwedge H$. Furthermore, $C(H)$ has a natural intersection form
induced by $<,>$, also denoted by $<,>$, and possesses the
involution
$$ t:C(H)\rightarrow C(H),$$
$$h_1\ldots h_k\mapsto h_k\ldots h_1,\,h_i\in H.$$
For $v\in C(H)$, $t(v)$ is the adjoint of the multiplication on the
left or on the right by $v$.
 Consider now $H_\mathbb{C}$, with its
Hodge decomposition
$$H_\mathbb{C}=H^{2,0}\oplus H^{1,1}\oplus H^{0,1}.$$
By assumption, the rank of $H^{2,0}$ is $1$, and equivalently the
dimension over $\mathbb{R}$ of the real vector space
$$(H^{2,0}\oplus H^{0,2})\cap H_\mathbb{R}$$
is $2$. By the second Hodge-Riemann bilinear relations, the
intersection form $<,>$ is positive definite on this $2$-plane.
Furthermore this $2$-plane is canonically oriented because it is
canonically isomorphic to the complex line $H^{2,0}$ via the map
$Re$. Choose an oriented orthonormal basis $e_1,\,e_2$ of this
$2$-plane. Then
$$e:=e_2e_1\in C(H_\mathbb{R})$$
does not depend on the choice of the basis. Furthermore we have
$$e^2=e_2e_1e_2e_1,$$
with
$$e_1e_2+e_2e_1=0,\,e_1^2=-1,\,e_2^2=-1.$$
Thus $e^2=-1$ and  left multiplication by $e$ defines a complex
structure on $C(H_\mathbb{R})$. This provides a complex torus
$$K(H)=C(H_\mathbb{R})/C(H),$$
which is the Kuga-Satake variety of $H$.

We want now to present this construction from our point of view of
weight $2$ Hodge structures endowed with a compatible  ring
structure.

We noticed already that as vector spaces, we have a canonical
identification
$$C(H_\mathbb{C})=\bigwedge H_\mathbb{C}.$$
The right hand side is a direct sum
$$\bigoplus_{k}\bigwedge^kH_\mathbb{C}.$$
On $\bigwedge^kH_\mathbb{C}$, there is a weight $2k$ Hodge
decomposition induced by the weight $2$ Hodge decomposition on
$H_\mathbb{C}$, with $H^{p,q}$ term given by
$$\bigoplus_{(r,s,t),2r+s=p,2t+s=q}\bigwedge^rH^{2,0}\otimes\bigwedge^sH^{1,1}\otimes
\bigwedge^tH^{0,2}.$$ As $rk\,H^{2,0}=1$, this decomposition has in
fact only three  terms, according to the value $0$ or $1$ given to
$r$ and $s$,  noticing that for $(r,s)=(1,1)$ and $(r,s)=(0,0)$, we are
in $H^{k,k}$. In other words, the induced Hodge structure on
$\bigwedge^kH$ has Hodge level $2$, that is can be shifted to a
weight $2$ Hodge structure, by shifting the bigrading by
$(1-k,1-k)$. Thus $H^{2,0}\wedge \bigwedge^{k-1}H^{1,1}$ is
assigned type $(2,0)$, $H^{0,2}\wedge \bigwedge^{k-1}H^{1,1}$ is
assigned type $(0,2)$, while $\bigwedge^{k}H^{1,1}$ and
$H^{2,0}\wedge H^{0,2}\wedge\bigwedge^{k-2}H^{1,1}$ are
assigned type $(1,1)$.

Taking the direct sum over  integers $k$, we get a weight $2$
Hodge structure on $C(H)=\bigwedge H$.

\begin{lemm} \label{hodgemult}This Hodge structure is compatible with the product on $C(H)$,
that is satisfies condition \ref{item1} of Theorem
\ref{trivialcenter}. \end{lemm} {\bf Proof.}  We have the Clifford
multiplication
$$ H\otimes C(H)\rightarrow C(H).$$ If we show that this map is a morphism of Hodge
structures of bidegree $(-1,-1)$, then the same will be true by
iteration (because $H$ generates $C(H)$) for the multiplication map
$$C(H)\otimes C(H)\rightarrow C(H).$$
But it is well-known (see \cite{lawsonmichel}, p. 25) that the Clifford
multiplication by $h\in H$ acting on $C(H)\cong\bigwedge H$
identifies to
$$ h\wedge-h\lrcorner ,$$
where $h\lrcorner $ acts on $\bigwedge H$ via the element of $H^*$
given by $<h,\cdot>$.

Now, let $h\in H^{2,0}=\eta\mathbb{C}\subset H_\mathbb{C}$. Then
$h\wedge $ annihilates the $(2,0)$-part of $\bigwedge^kH_\mathbb{C}$
which is equal to $\eta\wedge\bigwedge^{k-1}H^{1,1}$ and
$h\lrcorner$ also annihilates the $(2,0)$-part of
$\bigwedge^kH_\mathbb{C}$ because $<H^{2,0},F^1H_\mathbb{C}>=0$.

Next, the $(1,1)$-part of $\bigwedge^kH_\mathbb{C}$ is equal
to $$\eta\wedge\overline{\eta}\wedge\bigwedge^{k-2}H^{1,1}\oplus
\bigwedge^{k}H^{1,1}.$$ The map $h\wedge$ annihilates the first
term, and sends the second one in
$$\eta\wedge\bigwedge^{k-2}H^{1,1}$$
that is, in the $(2,0)$-part of $C(H)$. Furthermore, the map
$h\lrcorner$ annihilates the second
term, and sends the first one in
$$\eta\wedge\bigwedge^{k-2}H^{1,1}$$
that is, in the $(2,0)$-part of $C(H)$.

Finally, the $(0,2)$-part of
$\bigwedge^kH_\mathbb{C}$  is equal to
$\overline{\eta}\wedge\bigwedge^{k-1}H^{1,1}$, and $h\wedge$ sends
it to $$\eta\wedge\overline{\eta}\wedge\bigwedge^{k-1}H^{1,1}$$
which is contained in the $(1,1)$-part of $C(H)$, while $h\lrcorner$
sends it to $\bigwedge^{k-1}H^{1,1}$ which is contained in the
$(1,1)$-part of $C(H)$.

In other words, we proved that  Clifford multiplication by $h\in
H^{2,0}$ shifts the Hodge decomposition on $C(H_\mathbb{C})$ by
$(1,-1)$. One shows similarly that  Clifford multiplication by
$h\in H^{1,1}$ preserves the Hodge decomposition on $C(H)$, which
proves the claim.
 \cqfd
 It turns out that the assumption on the center used in the previous section  is not always satisfied by
  the Clifford algebra.
 However, it is quite easy to see directly in this case  that
 our definition of an associated weight $1$ Hodge structure on $C(H)$ still works in this case, that
 is, the factor $M$ in the decomposition (\ref{orthog}) is $0$. Thus we have
 the weight $1$ Hodge decomposition
 $$C(H_\mathbb{C})=W\oplus\overline{W}$$
 as in the previous section.

 To conclude, we show the following:
 \begin{prop} The Kuga-Satake construction of
 a complex structure on
 $C(H_\mathbb{R})$ coincide with our construction
 of a weight $1$ decomposition on $C(H_\mathbb{C})$
 (or equivalently a complex structure on
 $C(H_\mathbb{R})$) above.

 \end{prop}
 {\bf Proof.}
 Recall that our weight $1$ decomposition of $C(H_\mathbb{C})$ is
 given by
 given by
 $$C(H_\mathbb{C})=W\oplus\overline{W},$$
with $W=C(H)^{2,0}C(H_\mathbb{C})$. Observe now that, still denoting
by $\eta$ a generator of $H^{2,0}$, we have
$$C(H)^{2,0}C(H_\mathbb{C})=\eta C(H_\mathbb{C}).$$
Indeed, as $\eta\in C(H)^{2,0}$, the inclusion $\supset$ is clear.
The reverse inclusion comes from
$$C(H)^{2,0}=\eta\wedge\bigwedge H^{1,1}=\eta \cdot\bigwedge H^{1,1},$$
where on the right, the $\cdot$ stands for  Clifford multiplication
rather than exterior multiplication.

Next, consider the Kuga-Satake construction: the complex structure
$I$ on $C(H_\mathbb{R})$ is given here by multiplication by
$e=-e_1e_2$, where $e_1,\,e_2$ is an oriented orthonormal basis of
the real part of $H^{2,0}\oplus H^{0,2}$. Choosing $\eta$ in such a
way that $<\eta,\overline{\eta}>=2$, we may assume (because
$<\eta,\eta>=0=<e_1,e_1>-<e_2,e_2>+2i<e_1,e_2>$) that
$$e_1=Re\,\eta,\,e_2=Im\,\eta.$$
Furthermore, the weight $1$ decomposition on $C(H_\mathbb{C})$
associated  to the Kuga-Satake complex structure is determined by the complex subspace
$W'\subset C(H_\mathbb{C})$
$$C(H_\mathbb{C})=W'\oplus\overline{W'},$$
where $W'$ is by definition
the $i$-eigenspace of the complex structure operator $I$ in
$C(H_\mathbb{C})$, that is the subspace generated by the
$$w-iI(w),\,w\in C(H_\mathbb{C}).$$
As $I$ is Clifford multiplication on the left by $e$, $W'$ is also the subspace
generated by the
$$ (1-ie)w,\,w\in C(H_\mathbb{C}).$$
On the other hand, we have
$$\eta=e_1+ie_2$$
and thus
$$\eta\overline{\eta}=(e_1+ie_2)(e_1-ie_2)=-2-2ie_1e_2=-2(1-ie).$$
Hence we conclude that
$$W=Im\,\eta\subset
Im\,\eta\overline{\eta}=Im\,(1-ie)=W'.$$ By the equality of
dimensions, we now conclude that we have equality.
 \cqfd

In conclusion, we have split the Kuga-Satake construction into two
parts:

i)  The observation (studied from the point of view of its
cycle-theoretic implications in \cite{voisinsym}) that for a
$K3$-type Hodge structure, the induced Hodge structures on the
exterior powers of $H$ have level $2$.

ii) The construction of Theorem \ref{trivialcenter}, which works for
much more general compatible ring structures on polarized weight $2$
Hodge structures.
\begin{rema} {\rm  Lemma \ref{hodgemult} is also true for the exterior algebra
structure on $\bigwedge H$. From our point of view, the key reason
for which we need the Clifford algebra, is the existence of the
adjunction map $t$, which is not satisfied by the exterior algebra.
}
\end{rema}
\section{The general case \label{sec3}}
We pass now to rational coefficients. We consider as in the previous section a
polarized rational Hodge structure $H$ of weight $2$, which is also a
$\mathbb{Q}$-algebra, such that the product
$$H\otimes H\rightarrow H$$
is a morphism of Hodge structures of bidegree $(-1,-1)$, and such that there exists
an adjunction map
$$t:H\rightarrow H$$
such that (left or right) multiplication by  $t(h)$ is
adjoint to (left or right) multiplication by $h$, with respect to the polarization form
$<,>$.
Our first goal is to study the possible centers of such an algebra.
\begin{lemm}\label{trivialhodge} The center $K$ is a trivial sub-Hodge structure of $H$.
\end{lemm}
{\bf Proof.} It is obvious that $K$ is a sub-Hodge structure of $H$, because it is defined
as the kernel of the map
$$H\rightarrow Hom\,(H,H),$$
$$h\mapsto (a\mapsto ha-ah),$$
and this map is a morphism of Hodge structures of bidegree $(-1,-1)$.

To see that it is a trivial Hodge structure, let
$0\not=\alpha\in K^{2,0}\subset K\otimes \mathbb{C}$. Then
we have
$$u:=\alpha t(\overline{\alpha})\in H^{1,1},$$
and $u\not=0$ because $<u,1>=<\alpha,\overline{\alpha}> >0$ by the Hodge-Riemann bilinear relations.
Thus we have, again by the Hodge-Riemann bilinear relations:
$$ <u,\overline{u}>=<ut(\overline{u}),1>\not=0.$$
Hence $$ut(\overline{u})=\alpha t(\overline{\alpha})\alpha t(\overline{\alpha})\not=0.$$
But as $\alpha$ is central, this is $0$ because $\alpha^2=0$, which is a contradiction.
\cqfd
\begin{lemm}\label{fields} The center $K$ of $H$ is a product of number fields. The $t$-invariant part
$K^+$ of $K$ is a product of totally real number fields $K_i$. If $K_i^+\not=K_i$, then
$K_i$ is a field which becomes isomorphic to
$\mathbb{C}$ under any embedding of $K_i^+$ into $\mathbb{R}$.
\end{lemm}
{\bf Proof.} $K$ is a commutative $\mathbb{Q}$-algebra, which is by the previous Lemma
contained in $H^{1,1}_{\mathbb{R}}$. Thus the pairing $<,>$ restricts to
a negative definite pairing on $K$. Note also that $K$ is clearly invariant under $t$, as
$t$ is a morphism of $\mathbb{Q}$-algebras $H\rightarrow H^{op}$.
On the other hand
we know that
this pairing
is of the form
$$ <\alpha,\beta>=<\alpha t(\beta),1>.$$
We want to show that $K$ does not contain nilpotent elements.
If $\alpha^n=0$, then $(\alpha t(\alpha))^n=0$, and
 if $\alpha t(\alpha)=0$, then also $\alpha=0$, because $<\alpha,\alpha>=<\alpha t(\alpha),1>$.

Thus it suffices to show that $\alpha^n=0$ implies $\alpha=0$,  when $t(\alpha)=\alpha$.
But if $\alpha^n=0$, with $\alpha=t(\alpha)$ and $n=2m$, then
also $<\alpha^m,\alpha^m>=0$ which implies that $\alpha^m=0$. Thus the order
$n$ of $\alpha$ cannot be even, and if it is odd, $n=2m+1$, we have $n\leq m+1$, that is $m=0$.

Thus $K$ is a product of number fields. Note that the intersection form on
$K$ must be of the form
$$<\alpha,\beta>=tr_{K/\mathbb{Q}}(y\alpha t(\beta))$$
for some $y\in K$. (Here $y\in K$ is defined by the condition that
$$ <\alpha,1>=Tr_{K/\mathbb{Q}}(y\alpha),\,\forall \alpha\in K.)$$
Thus, as $<\alpha,1>=<t(\alpha),1>$, we must have $y\in K^+$.

We  show now that $K^+$ is a product of totally real
number fields. But we know
that the intersection form
$<\alpha\beta,1>$ on $K^+$ is negative definite, and that it is
of the form
$$<\alpha,\beta>=2tr_{K^+/\mathbb{Q}}(y\alpha\beta)$$
for some $y\in K^+$.

The algebra $K^+\otimes\mathbb{R}$ splits as a product of
quadratic extensions of $\mathbb{R}$, and we want to show that none
of these extensions can be $\mathbb{C}$. But
for any non-zero $y\in \mathbb{C}$, the quadratic form
$$<\alpha,\beta>=Tr_{\mathbb{C}/\mathbb{R}}(y\alpha\beta)$$
has signature $(1,1)$, which contradicts the fact that $<,>$ should be negative definite
on any factor of $K^+\otimes\mathbb{R}$.

To conclude, consider a component $K_i^+$ of $K^+$, given by
an idempotent $e_i\in K^+$. Let
$K_i:=e_i K$. Choose an imbedding $\sigma$
of $K_i^+$ into $\mathbb{R}$. As $K_i$ is a quadratic extension of
$K_i^+$, $K_i\otimes_{\sigma(K_i^+)}\mathbb{R}$ is a quadratic extension
$E$ of
$\mathbb{R}$, which is contained in $K\otimes_{\mathbb{Q}}\mathbb{R}\subset H^{1,1}_{\mathbb{R}}$.
Let this extension be given by $X^2=\lambda$.
The involution $t$ of this quadratic extension generates its  Galois group, and we have
the condition that
the quadratic form
$$<\alpha,\beta>=Tr_{E/\mathbb{R}}(y\alpha t(\beta))$$
for some $y\in \mathbb{R}$, is definite negative. But as $t(X)=-X$,
the matrix of this quadratic form in the base $(1,X)$ is
$$\begin{pmatrix}y&0\\
0&-y\lambda
\end{pmatrix}
$$
and thus we must have $\lambda<0$.

\cqfd

From now on, we assume that $K$ is a number field. This is possible because
the idempotents $e_i$ which give the decomposition of
$K$ into a product of number fields
$K_i$ are $t$-invariant Hodge classes in $H$ by Lemmas \ref{trivialhodge}
and \ref{fields}. Thus replacing
$H$ by $e_iH$, we still have a polarized Hodge structure, an adjunction map
$t$, and the compatibility of the product with the Hodge decomposition.

We define as in the previous section
$$W:=H^{2,0}H.$$
Then we have the decomposition (\ref{orthog})
$$ H_\mathbb{C}=(W\oplus\overline{W})\oplus M,$$
which is a orthogonal decomposition and an algebra decomposition.
$M$ is defined over $\mathbb{R}$, $M=M_\mathbb{R}\otimes\mathbb{C}$.
Note also that $M\subset H^{1,1}$, because $H^{2,0}\subset
W,\,H^{0,2}\subset\overline{W}$ and $(H^{2,0}\oplus
H^{0,2})^\perp=H^{1,1}$. Furthermore, $M$ is stable under $t$ and
the decomposition above is given by a central idempotent $e\in
K^+\otimes\mathbb{R}$. Thus, as $K^+$ is totally real,
$M_\mathbb{R}$ must be a sum
$$M_\mathbb{R}=\oplus_{\sigma\in\Sigma_M} H_\sigma,\,H_\sigma:= H\otimes_{\sigma(K^+)}\mathbb{R},$$
where $\Sigma_M$ is a certain set of imbeddings of $K^+$ into
$\mathbb{R}$. (Here we see $H\otimes_{\sigma(K^+)}\mathbb{R}$ as the
sub-algebra of $H_{\mathbb{R}}$ defined as the image  of the
idempotent $e_\sigma$ of $K^+\otimes_{\mathbb{Q}}\mathbb{R}$ given
by $\sigma$.)

\begin{prop} If $M\not=0$, then the algebra $H$ is a simple central
$K$-algebra.
\end{prop}
{\bf Proof.} We want to show that if $\Sigma_M$ is not empty, then
$H$ has no non trivial two-sided ideal. As $\Sigma_M\not=\emptyset$,
it clearly suffices  to show that the algebra
$H_\sigma$ has no non trivial two-sided ideal for $\sigma\in \Sigma_M$.
But the $\sigma\in \Sigma_M$ are characterized by the fact that
$H_\sigma\subset M$, which implies
$$H_\sigma\subset H^{1,1}.$$
Note that each $H_\sigma$ is invariant under $t$. Furthermore the
$t$-invariant part of the center of $H_\sigma$ is equal to
$\mathbb{R}$.

Let thus $I\subset H_\sigma$ be a two-sided ideal.
$I\cap t(I)$ is also a two-sided ideal which is $t$-invariant.
Furthermore, if $I\cap t(I)=0$, then
$I=0$. Indeed, if $x\in I$,
$xt(x)\in I\cap t(I)$. Thus, if $I\cap t(I)=0$, $xt(x)=0$. On the other hand,
we have $<x,x>=<xt(x),1>$ which is then also $0$. But as $x$ is real
of type $(1,1)$, this implies that
$x=0$ by the second Hodge-Riemann bilinear relations.

Thus we may assume that $I$ is $t$-invariant.
The orthogonal complement $J:=I^\perp$ of
$I$ in $H_\sigma$  with respect to $<,>$ is then also a $t$-invariant two-sided ideal of
$H_\sigma$.

On the other hand, the intersection form $<,>$ restricted to $I$ is non degenerate, because
 for $ I\subset H_\sigma\subset H^{1,1}$, and $<,>$ is negative definite on $H^{1,1}$. Hence
 we get an orthogonal decomposition of $H_\sigma$ into the sum of two two-sided ideals,
 or equivalently a decomposition
 $$H_\sigma=I\oplus J,$$
 as the direct sum of two sub-algebras.
 But then, as $I$ and $J$ are $t$-invariant, the idempotents associated to $I$ and $J$ are central and $t$-invariant.
 As the center of $H_\sigma$ is $\mathbb{R}$, it follows that either $I$ or $J$ is $0$.
 \cqfd
 \begin{coro} The conclusion
 of Theorem \ref{trivialcenter}  holds without the hypothesis \ref{item3} on the center,
 unless possibly when
 $H$ has a direct summand which is a simple central algebra $H_i$ over a number field
 $K_i$, where the decomposition is also a Hodge structure decomposition, and an orthogonal
 decomposition.
\end{coro}
\begin{ex}
{\rm We give here an example of
a rational polarized weight $2$ Hodge structure,
which admits a ${\mathbb Q}$-algebra structure, satisfying
conditions \ref{item1} and \ref{item2} of Theorem \ref{trivialcenter},
with center a number field $K$, and for which we have both
$$M\not=0,\,W\not=0,$$
where $W$ and $M$ are defined as in section \ref{sec1}.

We start from a rank $2$ vector space $V$ over the number field
$K:={\mathbb Q}(\sqrt{2},i)$, that is, $K={\mathbb
Q}[x,y]/(x^2=2,y^2=-1)$. Let $K^+:={\mathbb Q}(x)={\mathbb
Q}(\sqrt{2})\subset K$. On $K$, there is a unique $K^+$-bilinear
skew-symmetric form $\Omega$, which satisfies the property that
$$\Omega(1,y)=1.$$
On $V=K\oplus K$, consider the rational skew-symmetric bilinear form
:
$$\omega :=tr_{K+/{\mathbb Q}}\Omega_1+tr_{K+/{\mathbb Q}}x\Omega_2,$$
where $\Omega_i$ is $\Omega$ on the $i$-th factor.

We want now to put a weight $1$ Hodge structure on
$V$, polarized by $\omega$ and admitting $K$ as an
endomorphism algebra.
This is done as follows:

The space $V\otimes_{\mathbb Q}{\mathbb R}$ splits as the
direct sum
$$V_{\sqrt{2}}\oplus V_{-\sqrt{2}}$$
corresponding to the two embeddings of $K^+$ in
${\mathbb R}$ sending $x$ to
$\sqrt{2}$ or to $-\sqrt{2}$ respectively. Each
term of this decomposition admits the action of $y$, with $y^2=-1$.

We want  to put a complex structure on $V_{\mathbb R}$, given by an
operator of complex structure $I$,  which leaves invariant $\omega$,
and is such that $\omega(x,Ix)>0$ for $0\not= x\in V_{\mathbb R}$.
We define for this $I=y$ on $V_{\sqrt{2}}$, while on the factor
$V_{-\sqrt{2}}$, the operator $I$ is defined as follows:

The space
$$K_{\mathbb R}:=K\otimes{\mathbb R}$$
is equal to the quadratic extension of $K^+_{\mathbb
R}:=K^+\otimes{\mathbb R}$ given by $y^2=-1$. The algebra
$K^+_{\mathbb{R}}$ splits  as a sum of two copies of ${\mathbb R}$,
$K^+_{\sqrt{2}}$ and $K^+_{-\sqrt{2}}$, where $x$ acts by
multiplication by $\sqrt{2}$ and $-\sqrt{2}$ respectively. The first
copy is canonically isomorphic to ${\mathbb R}$ generated by the
idempotent $1^+:=\frac{x+\sqrt{2}}{2\sqrt{2}}$ and the second copy
is similarly generated by the idempotent
$1^-:=\frac{x-\sqrt{2}}{-2\sqrt{2}}$.

This makes (via the choice $y=i$) the algebra $K_{\mathbb R}$
canonically isomorphic to the sum of two copies $K_{\sqrt{2}}$,
$K_{-\sqrt{2}}$ of ${\mathbb C}$, and thus $V_{\mathbb R}\cong
K_{\mathbb R}^2$ is a direct sum of four copies of ${\mathbb C}$. It
is not hard to compute that the extended form $\omega$ in the
canonical coordinates $z_1,\ldots,z_4$ on ${\mathbb C}^4$ is equal
to
\begin{eqnarray}\label{omega}\frac{i}{2}(dz_1\wedge d\overline{z_1}+ dz_2\wedge
d\overline{z}_2+\sqrt{2}dz_3\wedge d\overline{z}_3
-\sqrt{2}dz_4\wedge d\overline{z}_4). \end{eqnarray} Thus it is of
type $(1,1)$ for the complex structure given by $y$, but not
positive.

We will define our $I$ to be equal to $i$ on the first three
factors and to $-i$ on the last factor. This has the effect of
exchanging the coordinates $z_4$ and $\overline{z}_4$, and then by
formula (\ref{omega}),  $\omega$ becomes positive of type $(1,1)$
with respect to $I$. By construction, $I$ commutes with $x$ and $y$, as it preserves their eigenspace
decomposition.

Having this, we will now consider as polarized Hodge structure
$$H:=End_K(V).$$
This is a sub-Hodge structure of
$End\,(V)$, because $K$ acts as endomorphisms of
the weight $1$ Hodge structure on $V$.

This is stable under adjunction with respect to $\omega$, because
the action of $K$ on $V$ is such that multiplication by $x$ is
self-adjoint with respect to $\omega$, while multiplication by $y$
is anti-self-adjoint.

On the other hand,
consider the decomposition
$$H_{\mathbb R}=H_{\sqrt{2}}\oplus H_{-\sqrt{2}}.$$
Then $H_{\sqrt{2}}= End_{\mathbb C}(V_{\sqrt{2}})$ and
$H_{-\sqrt{2}}=End_{\mathbb C}(V_{-\sqrt{2}})$, where the
${\mathbb C}$-structures are given by the action of
$iy$ on both terms.
As $I=y$ on $V_{\sqrt{2}}$, the term
$H_{\sqrt{2}}$ is made of endomorphisms commuting with $I$, thus of type
$(1,1)$, while as $I\not=y$ on $V_{-\sqrt{2}}$,
the term $H_{-\sqrt{2}}$ is not of type $(1,1)$.
Thus in this case both $M=H_{\sqrt{2}}$ and
$W\oplus\overline{W}=H_{-\sqrt{2}}$ are non zero.

}

\end{ex}

We now come back to the general case of a rational Hodge structure
$H$, polarized by an intersection form $<,>$, and endowed with a
${\mathbb Q}$-algebra structure, satisfying conditions \ref{item1}
and \ref{item2} of Theorem \ref{trivialcenter}. We also assume that
the $t$-invariant part of the center  is a number field $K^+$.

Our goal is to show the following:
\begin{theo}\label{nontrivialcenter}
There exists a polarized weight $1$ Hodge on
$H$, such that the multiplication on the left
$$H\rightarrow Hom\,(H,H)$$
is a morphism of weight $2$ Hodge structures,
where the weight $2$ Hoge structure on the
left is the initial one, and the weight $2$ Hodge structure
on the right is induced by the weight $1$ Hodge structure
on $H$.
\end{theo}

\begin{rema}{\rm      The main defect of this construction is the fact
that it is not unique, and does not satisfy the universal property
of Proposition \ref{universal}.}
\end{rema}
{\bf Proof.} Recall that
we have the decomposition
$$H_{\mathbb R}=(W\oplus\overline{W})_{\mathbb R}\oplus M_\mathbb{R},$$
which is a orthogonal decomposition, an
algebra decomposition, and is compatible
with Hodge decomposition.
$M$ is a sum of
factors
$$M=\oplus_{\sigma\in \Sigma_M} H_\sigma,$$
where
$H_\sigma:=H\otimes_{\sigma(K^+)}{\mathbb R}$
is a subalgebra of $H_{\mathbb R}$.
There is already a complex structure on the first term,
given by the isomorphism of real vector spaces
$$Re:W\cong(W\oplus\overline{W})_{\mathbb R},$$
and we simply have to put a complex structure $I_\sigma$ on each
component $H_\sigma$ for $\sigma\in \Sigma_M$.

Note that we want the left  multiplication map
$$H\rightarrow Hom\,(H,H)$$ to be a morphism of Hodge structures.
This implies that
$$H_\sigma\rightarrow Hom\,(H_\sigma,H_\sigma)$$ has to
be a morphism of real Hodge structures. But, for $\sigma\in
\Sigma_M$, we know that $H_\sigma$ is of type $(1,1)$. Thus
multiplication on the left on $H_\sigma$ by any $h\in H_\sigma$ has
to be of type $(1,1)$, which means that it commutes with the complex
structure operator $I_\sigma$. But as our algebra has an unit, this
implies in turn that $I_\sigma$ has to be the multiplication on the
right by some element $m_\sigma\in H_\sigma$, satisfying the
condition that $m_\sigma^2=-1$.

Furthermore, we want that our weight $1$ Hodge structure
is polarized, with a polarization of the form
$$\omega_a(x,y)=<x,ya>,$$
for some $a\in H$ satisfying $t(a)=-a$. (This was the form
chosen for the polarization on the $W$-term, and
since the polarization must be  rational, we do not have another choice
here.)

The first condition for $\omega_a$ to polarize the real Hodge
structure on $H_\sigma$ is the fact that $\omega_a$ is of type
$(1,1)$ for the complex structure $I_\sigma$. Equivalently, for any
$x,\,y\in H_\sigma$,
$$<x,ya>=<I_\sigma(x),I_\sigma(y)  a>=<xm_\sigma,ym_\sigma a>=
<x,ym_\sigma at(m_\sigma)>.$$

This implies that
\begin{eqnarray} \label{relation} m_\sigma at(m_\sigma)= a.
\end{eqnarray}

Let us now distinguish the cases where $K=K^+$, $K\not=K^+$.

a) {\it Case $K\not= K^+$.} In this case, we proved in Lemma \ref{fields}
that $K_\sigma\cong \mathbb{C}$. Choosing such an isomorphism gives
$i_\sigma\in K_\sigma$, with $i_\sigma^2=-1$. Thus we  have  the operator of
complex structure acting on $H_\sigma$  by multiplication by
$i_\sigma$.

As $i_\sigma$ is in the center of $H_\sigma$, and satisfies
$t(i_\sigma)=-i_\sigma$, the relation (\ref{relation}) for
$m_\sigma=i_\sigma$ is certainly satisfied for all $a$.

It remains to see that with this operator of complex structure
$I_\sigma$, the corresponding real weight $1$ Hodge structure on
$H_\sigma$  is polarized by $\omega_a$ for an adequate $a\in H$
satisfying $t(a)=-a$. As we know already that $\omega_a$ is of type
$(1,1)$ for $I_\sigma$, we have only to verify that
$$\omega_a(x,I_\sigma x)>0,$$
for $0\not=x\in H_\sigma$. By definition, this is equal to
$$<x,xi_\sigma a>.$$

For $a=i_\sigma\in H_\sigma\subset H_{\mathbb R}$, this is equal to
$-<x,x>$ which is positive by the second Hodge-Riemann bilinear
relations, because $H_\sigma\subset H^{1,1}_{\mathbb R}$.

Thus it follows that it remains positive for any $a$ in a
neighbourhood of $i_\sigma$ in $H_{\mathbb R}$ satisfying $t(a)=-a$,
and in particular for a rational such $a\in H\subset H_{\mathbb R}$.

b){\it Case $K=K^+$.}
Here the $a$ will be a  fixed rational element of $H$ such that $\omega_a$
satisfies the positivity conditions on the $W$-components (the
existence of which was shown in the previous section), and on the
components $H_\sigma$ of the previous type a).
Note that the multiplication (on the right or on the left) by $a$ on $H$
is an isomorphism, because it is $K$-linear and for any
imbedding $\tau:K\hookrightarrow {\mathbb R}$  such that $\tau\not\in
\Sigma_M$, it induces an automorphism of $H_\tau$, because $\omega_a$ is non-degenerate
on $H_\tau$.

(The existence of imbeddings $\tau\not\in \Sigma_M$ follows from the fact
that
$W\not=0$ by assumption.)

We consider the commutative $K$-subalgebra of
$H$ generated by $a$:
$$K_a:=K[a]\subset H.$$
This subalgebra is invariant under $t$, because $t(a)=-a$.
Next,  the
subalgebra  $K_{a,\sigma}:=K_a\otimes_{\sigma(K)}{\mathbb R}$
of $H_\sigma$ is contained in $H^{1,1}$, and thus satisfies the property
that the intersection form
$$<x,y>=<xt(y),1>$$
is negative definite on $K_{a,\sigma}$.

This implies as in the proof of
Proposition \ref{fields} that the $t$-invariant part $K_{a,\sigma}^+$
of $K_{a,\sigma}$ is a sum of copies
${\mathbb R}_\rho$  of ${\mathbb R}$,
and that the corresponding decomposition
of $K_{a,\sigma}$ (given by the action of the idempotents of
$K_{a,\sigma}^+$) is a decomposition as a sum of copies
${\mathbb C}_\rho$  of
${\mathbb C}$, where $t$ acts as complex conjugation
on each ${\mathbb C}_\rho$. Furthermore, as
$t(a)=-a$ and multiplication by $a$ is an isomorphism on $K_a$,
each $a_\rho$ can be written uniquely  as
$\lambda_\rho i_\rho$, where $\lambda_\rho$ is a positive real number,
and $i_\rho\in {\mathbb C}_\rho$ satisfies
 $i_\rho^2=-1_\rho$.

Let us define $m_\sigma$ by
$$m_\sigma=\sum_\rho i_\rho.$$
We have $m_\sigma^2=-\sum_\rho1_\rho=-1$. Furthermore
$t(m_\sigma)=-m_\sigma$, as $t(i_\rho)=-i_\rho$, for all $\rho$.

Finally, as $m_\sigma\in K_{a,\sigma}$,
one has
$m_\sigma a=a m_\sigma$, and combining these three facts, we conclude that
(\ref{relation}) is satisfied.

Thus the multiplication on the right by $m_\sigma$ defines a complex
structure
$I_\sigma$  on $H_\sigma$, which satisfies the property that $\omega_a$
is of type $(1,1)$ for $I_\sigma$.
In order to conclude that we have a real weight one
polarized Hodge structure on $H_\sigma$, we have to check
the positivity property
$$\omega_a(x,I_\sigma(x))>0,\,\forall 0\not=x\in H_\sigma.$$
But as the numbers $\lambda_\rho $ are positive,
one has
$$a_\rho=n_\rho^2 m_\rho,$$
for some $n_\rho\in {\mathbb R}$.
Thus, letting
$$n:=\sum_\rho n_\rho,$$
we have
$$a_\sigma=n^2 m_\sigma,$$
with $n\in K_{a,\sigma}$,  $t(n)=n$.
But then we have
$$\omega_a(x,I_\sigma(x))=
<x,I_\sigma(x) a>=<x,xm_\sigma a>$$
$$=<x,x n^2
m_\sigma^2>=-<nx,nx>, $$
where the last inequality holds because
$m_\sigma^2=-1$ and $t(n)=n$.
As $n$ is non-degenerate, the second Hodge bilinear relations
show that $<nx,nx>$ is negative for all $0\not=x\in H_\sigma$,
which is what we wanted.

 \cqfd
\begin{rema}{\rm Assume that either  $\Sigma_M$ is empty, or
that the center $K$ of $H$ satisfies $K\not=K^+$. In the first case,
we constructed the weight $1$ Hodge structure by defining its
$H^{1,0}$-part to be $H^{2,0}H_{\mathbb C}$, which is a right ideal.
Thus multiplication on the right by elements of $H$ are morphisms of
weight $1$ Hodge structure.

In the second case, the same is true if we choose for complex
structure operator on the $H_\sigma,\,\sigma\in H_\sigma$, the
multiplication by $i_\sigma\in K\otimes_{\sigma(K^+)}{\mathbb R}$.
Indeed, as $i_\sigma$ is central, this multiplication commutes with
right multiplication with elements of $H$.

Thus, in both cases,   as in the Kuga-Satake case, we can construct
an abelian variety, which admits $H$ as a sub-Hodge structure of
weight $2$, and also $H^{op}$ as a ring of endomorphisms.}
\end{rema}

\end{document}